\documentclass{amsart}

\usepackage{graphicx}
\usepackage{amsfonts}

\begin{document}

\title[THE CORRECT CLASSIC GENERALIZED LEAST-SQUARES ESTIMATOR OF ...]
{THE CORRECT CLASSIC
GENERALIZED LEAST-SQUARES ESTIMATOR OF 
AN UNKNOWN CONSTANT MEAN OF RANDOM FIELD}

\author{T. SUS{\L}O}
\email{tomasz.suslo@gmail.com}

\begin{abstract} 
The aim of the paper is to derive for 
the negative correlation function with a time parameter
an asymptotic disjunction $\lim_{j \rightarrow \infty}\omega^i_j v_i$ of 
the numerical generalized least-squares estimator $\omega^i_j v_i$ of 
an unknown constant mean of 
random field in fact the correct classic generalized
least-squares estimator of an unknown constant mean of the field.
\end{abstract}

\maketitle
 
\thispagestyle{empty}

\setcounter{footnote}{0}
\renewcommand{\thefootnote}{\alph{footnote}}

\vspace*{4pt}
\normalsize\baselineskip=13pt  
\section{Introduction}
\noindent
The best linear unbiased generalized (estimation) statistics 
$\hat{V}_j=\sum_{i=1}^n\omega^i_j V_i=\omega^i_j V_i$  
of the random field $V_j;~j \subset i=1,\ldots,n$ at $j \ge n+1$
with unknown constant mean $m$ and variance $\sigma^2$
that fulfils 
the constraint
$$
\lim_{j \rightarrow \infty}
E\{[V_j-\omega^i_j V_i]^2\}
=
\sigma^2
=
E\{[V_j-m]^2\} 
$$
is the classic best linear unbiased generalized statistics 
for finite $n$ and $j \rightarrow \infty$ of 
an unknown constant mean $m=E\{V_j\}$ of the field $V_j$
with the classic generalized least-squares estimator 
$\lim_{j \rightarrow \infty}\omega^i_j v_i$ of 
an unknown constant mean of the field 
and with constrained variance of the best linear unbiased generalized (estimation) statistics 
$\lim_{j \rightarrow \infty}E\{[\omega^i_j V_i-m]^2\}$ of the field
as its variance (a mean squared error of mean estimation). \\
The best linear unbiased generalized (estimation) statistics that 
fulfils (on computer) 
the constraint 
$$
E\{[V_j-\omega^i_j V_i]^2\}
=
\sigma^2
=
E\{[V_j-m]^2\} 
$$
is the numerical best linear unbiased generalized statistics 
for finite $n$ at finite $j$ of 
an unknown constant mean $m=E\{V_j\}$ of 
the field $V_j$ with the numerical generalized least-squares 
estimator $\omega^i_j v_i$ of 
an unknown constant mean of the field 
and with constrained variance of the best linear unbiased generalized (estimation) statistics
$E\{[\omega^i_j V_i-m]^2\}$ of the field as its variance. \\
Since the classic best linear unbiased generalized statistics
for finite $n$ and $j \rightarrow \infty$ of
an unknown constant mean $m=E\{V_j\}$ of the field $V_j$
is an asymptotic disjunction for $j \rightarrow \infty$ of
the numerical best linear unbiased generalized statistics for
finite $n$ at finite $j$ of an unknown constant mean $m=E\{V_j\}$ of
the field $V_j$ 
then the correct classic generalized least-squares estimator
$\lim_{j \rightarrow \infty} \omega^i_j v_i$ of
an unknown constant mean $m$ of the field is an asymptotic
disjunction for $j \rightarrow \infty$
of the numerical generalized least-squares estimator
$\omega^i_j v_i$ of an unknowm constant mean $m$ of
the field.

\section{The correlation function and the estimators}
\noindent
From the so-called semi-variogram  $\gamma(h)$ 
for a stationary random field $V_j$
with an unknown constant mean $m=E\{V_j\}$,
variance $\sigma^2=E\{V_j^2\}-E^2\{V_j\}$
and covariance function $C(h)=E\{V_j V_{j+h}\}-E\{V_j\}E\{V_{j+h}\}$
\begin{eqnarray}
\gamma(h) 
&=&  
\frac{1}{2}E\{(V_j-V_{j+h})^2\} \nonumber \\  
&=&
\frac{1}{2}E\{V_j^2-2V_jV_{j+h}+V_{j+h}^2\} \nonumber \\
&=&
E\{V_j^2\}-E\{V_j V_{j+h}\} \nonumber \\
&=&
E\{V_j^2\}-E^2\{V_j\}-(E\{V_j V_{j+h}\}-E^2\{V_j\}) \nonumber \\
&=&
E\{V_j^2\}-E^2\{V_j\}-(E\{V_j V_{j+h}\}-E\{V_j\}E\{V_{j+h}\}) \nonumber \\
&=&
\sigma^2-C(h) \ge 0 \nonumber
\end{eqnarray}
we get the absolute value of the correlation function $\rho(h)$
$$
|\rho(h)| 
=
C(h)\slash C(0)=C(h)\slash\sigma^2=1-\gamma(h)\slash\sigma^2 \ .
$$
Since the correlation function is non-increasing
then the (first) experimental correlogram 
\begin{equation}
|\hat{\rho}(h)|=1-\hat{\gamma}(h)\slash\hat{\sigma}^2
\label{c1}
\end{equation}
should be computed for non-decreasing outcomes for $h \le d$ 
$$
\hat{\sigma}^2
=
\hat{\gamma}(d)\ge\hat{\gamma}(d-1)\ge\hat{\gamma}(d-2)\ge
\ldots\ge\hat{\gamma}(1) > \hat{\gamma}(0)=0 
$$
of the experimental semi-variogram for a time series
of the length $n$ 
$$
\hat{\gamma}(h)
=\frac{1}{2}\frac{1}{(n-h)}\sum_{j=1}^{n-h}(v_{j}-v_{j+h})^2
\quad
h=0,\ldots,n-1 \ .
$$ 
On the other hand from definition of the covariance function 
$$
C(h)=E\{V_j V_{j+h}\}-E\{V_j\}E\{V_{j+h}\}
$$
we get for a time series of the length $n$
$$
\hat{C}(h)
=
\frac{1}{n-h}\sum_{j=1}^{n-h}v_jv_{j+h}
-  
\frac{1}{(n-h)^2}
\sum_{j=1}^{n-h}v_j
\sum_{j=1}^{n-h}
v_{j+h} 
$$ 
and the (second) experimental correlogram 
\begin{equation}
|\hat{\rho}(h)|=\hat{C}(h)\slash \hat{C}(0)
\label{c2}
\end{equation}
for $h=0,\ldots,d$.

\section{The correct classic generalized least-squares estimator of 
an unknown constant mean of the field}
\noindent 
The attached source code ``combo.pas'' 
let us find (see~Tab.~\ref{Tab1}) 
for the negative correlation function
with the time parameter $t=n+1,\ldots,n+s$
\begin{equation}
\rho(\Delta_{ij})=\left\{
        \begin{array}{ll}
        -1 \cdot {\displaystyle t}^{\displaystyle -\Theta[\Delta_{ij}\slash t]^2},& 
        \qquad \mbox{for}~~h=\Delta_{ij}=|i-j|>0,\\
        +1, & \qquad \mbox{for}~~h=\Delta_{ij}=|i-j|= 0,\\
        \end{array}
        \right. 
\label{cf}
\end{equation}
where $\Theta$ is derived by Levenberg-Marquardt fit 
$$
|\rho(h)|
=
{\displaystyle n}^{\displaystyle -\Theta[h\slash n]^2} \ , 
$$ 
to the experimental correlograms~(\ref{c1}) and (\ref{c2}),   
the asymptotic disjunction $\lim_{j \rightarrow \infty} \omega^i_j v_i$ of 
the numerical generalized least-squares estimator 
$\omega^i_j v_i$ of an unknown constant mean $m$ of the field 
$V_j;~j \subset i=1,\ldots,n$ 
in fact the correct classic generalized least-squares estimator  
of an unknown constant mean $m$ of the field. \\ 
The variance of the best linear unbiased generalized (estimation) statistics
\begin{equation}
E\{[\omega^i_jV_i-m]^2\}
=
- \sigma^2 (\omega^i_j \rho_{ij}-\mu_j^{\it 1})
\label{mseofme}
\end{equation}
under the constraint
$$
E\{[V_j-\omega^i_j V_i]^2\}
=
\sigma^2
=
E\{[V_j-m]^2\} 
$$
given by the numerical approximation to the root of equation
\begin{equation}
\omega^i_j \rho_{ij}+ \mu_j^{\it 1}=0  \ ,
\label{constraint}
\end{equation}
where (kriging algorithm)
$$
\begin{array}{cccccl}
\underbrace{
\left[
\begin{array}{c}
\omega_j^{\it 1} \\
\vdots \\
\omega_j^{\it n} \\
\mu^{\it 1}_j \\
\end{array}
\right]
}_{(n+1) \times 1}
&
=
{\underbrace{
\left[
\begin{array}{cccc}
\rho_{\it 11} & \ldots & \rho_{\it 1n} & 1 \\
\vdots & \ddots & \vdots & \vdots \\
\rho_{\it n1} & \ldots & \rho_{\it nn} & 1 \\
1 & \ldots & 1 & 0 \\
\end{array}
\right]}_{(n+1) \times (n+1)}}^{-1}
&
\cdot
&
\underbrace{
\left[
\begin{array}{c}
\rho_{\it 1j} \\
\vdots \\
\rho_{\it nj} \\
1 \\
\end{array}
\right]
}_{(n+1) \times 1} \ ,
\end{array}
$$
becomes a variance of the best linear unbiased generalized statistics of 
an unknown constant mean $m$ of random field with 95$\%$ confidence intervals
for a mean estimation (from~(\ref{mseofme}))
$$
\hat{m}=\omega^i_j v_i \pm 1.96 \sqrt{-\hat{\sigma}^2
(\omega^i_j \rho_{ij}-\mu_j^{\it 1})} \ ,
$$
where
$$
\hat{\sigma}^2=\omega^i_j v^2_i-(\omega^i_j v_i)^2 \ ,
$$
if (\ref{constraint}) holds.
\begin{table}[hb!]
\centerline{
\begin{tabular}{|c|c|c|c|c|c|c|c|}
\hline
     &       &     &        &         &       &  \\
Fig. & Index & $n$ & $\Theta$ & $t=n+s$ & $j=u$ & 
$\omega^i_j v_i$ \\ 
\hline
1 & FTSE 100& $ 132 $ & $0.83283$ & 238 & 426 & 8463.42 \\
\hline
2 & DJI & $ 113 $ & $0.93670$ & 203 & 356 & 13603.87 \\
\hline
3 & S\&P 500& $ 102 $ & $0.93569$  & 174 & 270 & 1788.80 \\   
\hline            
\end{tabular}}
\caption{\label{Tab1}\small 
The numerical generalized least-squares estimator 
$\omega^i_j v_i$ of an unknown constant 
mean $m=E\{V_j\}$ of the field $V_j$ for the 
negative correlation function~(\ref{cf})
with the time parameter $t=n+1,\ldots,n+s$ as the correct classic 
generalized least-squares estimator 
of an unknown constant mean of the field
for final $t=n+s$ at final $j=u$.}
\end{table}
\begin{figure}[t!] 
\vspace*{13pt}
\includegraphics[width=12cm,height=6cm]{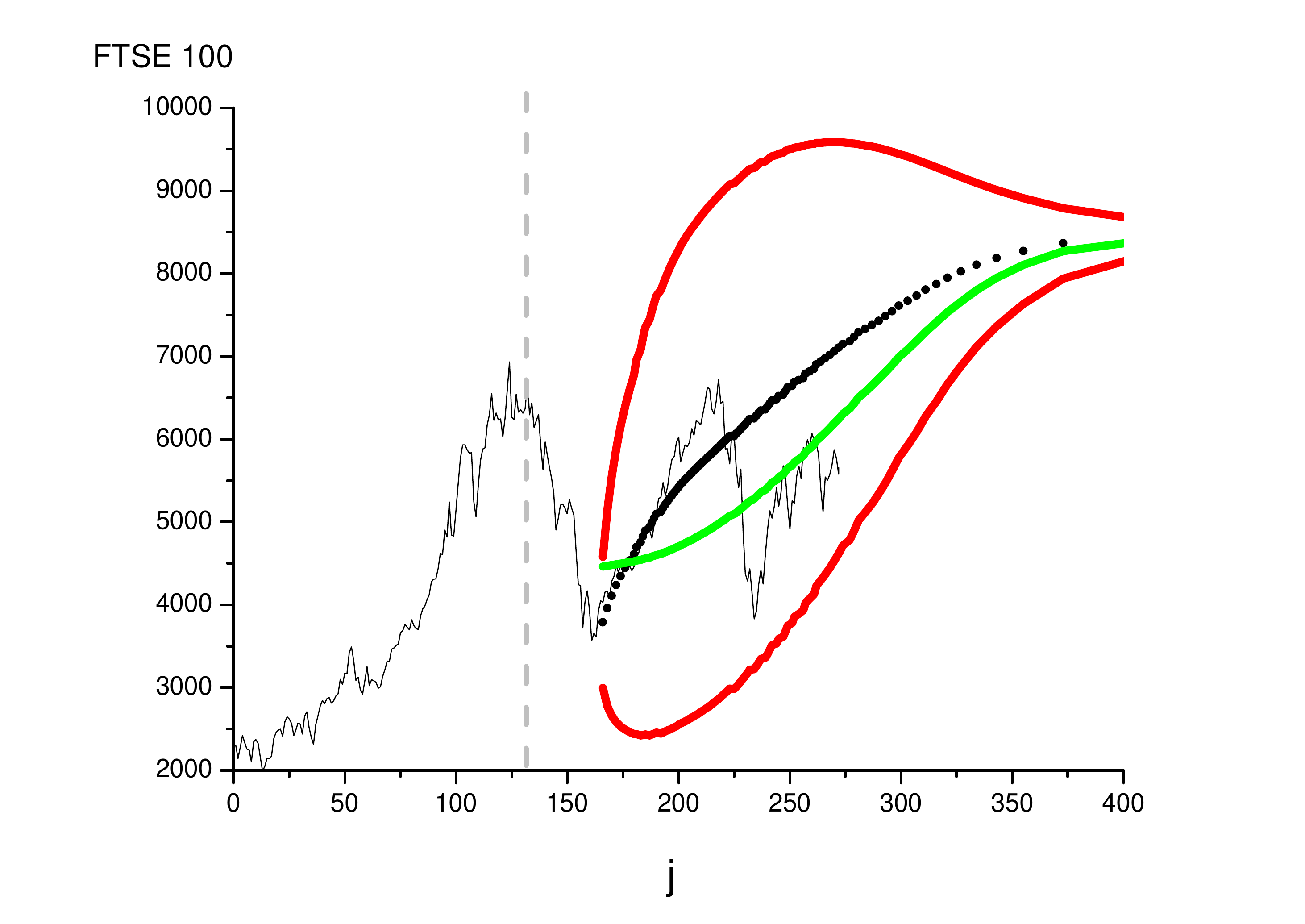}
\caption{\label{Fig1}\small 
FTSE 100 from 1 September 1989 up to 1 May 2012 (272 monthly close quotes). 
The numerical (black dots) generalized least-squares estimator 
$\omega^i_j v_i$ of 
an unknown constant mean of the field 
with 95\% confidence intervals (red lines) for a mean estimation is compared 
for the negative correlation function~(\ref{cf}) with the parameter  
$t=n+1,\ldots,n+s$ at finite $j \ge n+1$ to the classic (grey line) 
generalized least-squares estimator of an unknown constant mean 
$\lim_{j \rightarrow \infty}\omega^i_j v_i$ of the field. 
The asymptotic limit of 
the classic generalized least-squares 
statistics of an unknown constant mean of the field 
is fulfilled for final $t=238$ at final $j=426$. 
Dashed vertical line denotes $n=132$.
} 
\end{figure}
\begin{figure}[b!] 
\vspace*{13pt}
\includegraphics[width=12cm,height=6cm]{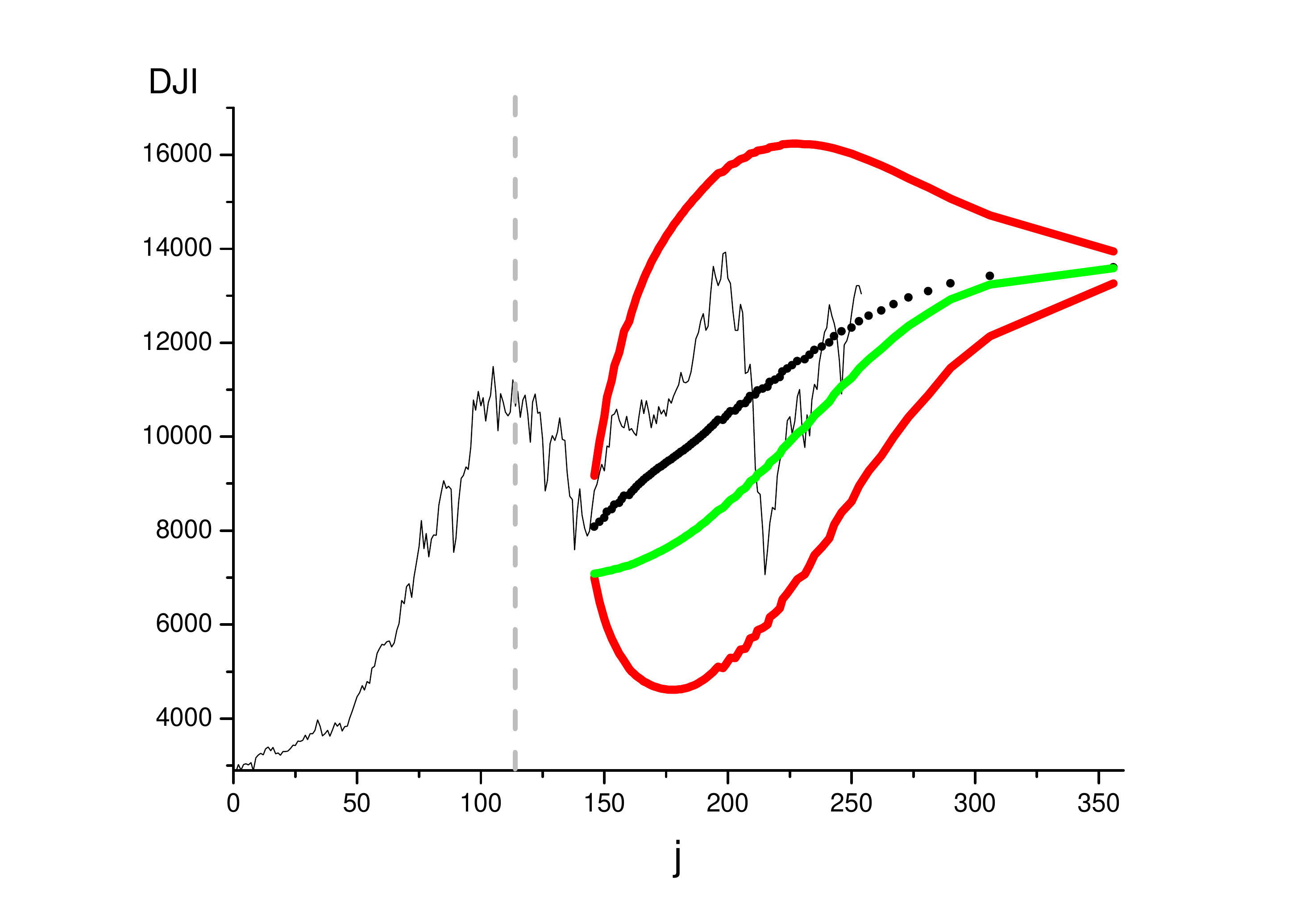}
\caption{\label{Fig2}\small 
DJI from 1 April 1991 up to 1 May 2012 (254   monthly close quotes). 
The same as in Fig.~\ref{Fig1}.
The asymptotic limit of 
the classic generalized least-squares 
statistics of an unknown constant mean of the field 
is fulfilled for final $t=203$ at final $j=356$. 
Dashed vertical line denotes $n=113$.
} 
\end{figure}
\begin{figure}[t!] 
\vspace*{13pt}
\includegraphics[width=12cm,height=6cm]{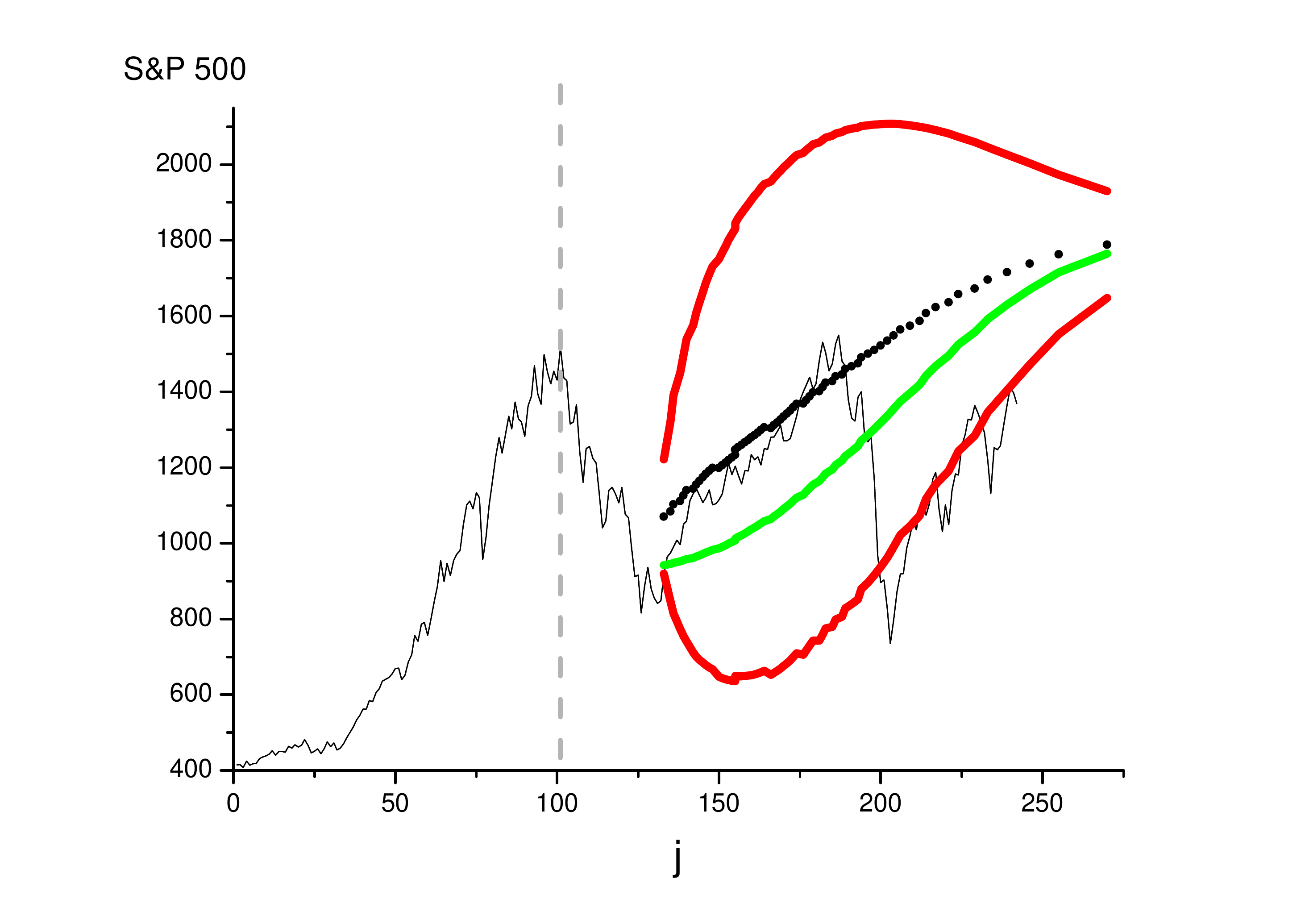}
\caption{\label{Fig3}\small 
S\&P 500 from 3 April 1992 up to 1 May 2012 (242 monthly close quotes). 
The same as in Fig.~\ref{Fig1}.
The asymptotic limit of 
the classic generalized least-squares 
statistics of an unknown constant mean of the field 
is fulfilled for final $t=174$ at final $j=270$. 
Dashed vertical line denotes $n=102$.
} 
\end{figure}

\end{document}